\documentclass[10pt,a4paper]{article}
\usepackage{amssymb}
\usepackage{}
\usepackage{colortbl}
\usepackage{amssymb,amsmath,graphicx,stmaryrd}
\usepackage{array,cases,dsfont,graphicx,texdraw}
\usepackage{graphics} 
\usepackage{epsfig} 
\usepackage{subfigure}

\newenvironment{Proof}{\noindent{\em Proof:\/}}{\hfill $\Box$\par}

\newtheorem{Remark}{Remark}[section]

\newtheorem{Definition}{Definition}[section]
\newtheorem{Theorem}{Theorem}[section]

\newtheorem{Lemma}{Lemma}[section]
\newtheorem{Property}{Property}[section]
\newtheorem{Assumption}{Assumption}[section]

\newcommand{\EQQ}{\begin{eqnarray*}}
\newcommand{\ENN}{\end{eqnarray*}}
\newcommand{\EQ}{\begin{eqnarray}}
\newcommand{\EN}{\end{eqnarray}}
\newcommand{\bd}{\begin{Definition}
\begin{rm} }
\newcommand{\ed}{ \end{rm}
\end{Definition} }

\begin{document}

\title{\Large {\bf Internal Model Approach to Cooperative Robust Output Regulation for Linear Uncertain Time-Delay Multi-Agent Systems}}

\author{Maobin Lu and Jie Huang
\thanks{This work has been supported in part by the Research Grants
Council of the Hong Kong Special Administration Region under grant
No. 412813, and in part by the National Natural Science Foundation
of China under grant No. 61174049.}
\thanks{M. Lu and J. Huang are with Shenzhen Research Institute, The Chinese University of Hong Kong, Shenzhen, China,  and Department of Mechanical and Automation Engineering, The Chinese University of Hong Kong, Shatin, N.T., Hong Kong. Email: mblu@mae.cuhk.edu.hk, jhuang@mae.cuhk.edu.hk.
        }%
}

\maketitle
\begin{abstract}

In this paper, we study the cooperative robust output regulation problem for linear uncertain multi-agent
systems with both communication delay and input delay by the distributed internal model approach. The problem includes the leader-following consensus problem of linear multi-agent systems with time-delay as a special case. We first generalize the internal model design method to systems with both communication delay and input delay. Then,
under a set of standard assumptions, we have obtained the solution of the problem via both the state feedback control and the output feedback control. In contrast with the existing results, our results apply to general uncertain linear multi-agent systems, accommodate a large class of leader signals, and  achieve the asymptotic tracking and disturbance rejection at the same time.

\end{abstract}

\section{INTRODUCTION}

In this paper, we consider the cooperative robust output regulation for
linear uncertain time-delay systems of the following form:
\begin{equation}\label{sys}
\begin{split}
\dot{x}_i(t) &= \bar{A}_i x_i(t) +  \bar{B}_i  u_i(t-\tau_{con}) + \bar{E}_i v(t), t \geq 0\\
y_i(t) &= \bar{C}_i x_i(t)  , t \geq 0\\
x_i(\theta) &=   x_{i0}(\theta), \     \theta \in [-\tau_{con},0]
\end{split}
\end{equation}
where $x_i(t) \in R^n$, $y_i(t) \in R^p$, and $u_i(t) \in R^m$ are the system state, measurement output, and control input of the $i^{th}$ subsystem, $\tau_{con}\geq 0$ is the input delay, and $v(t) \in R^q$ is the exogenous signal representing the reference
input to be tracked or/and disturbance to be rejected and is assumed to be
generated by the exosystem of the form
\begin{equation}\label{exo}
\begin{split}
\dot{v}(t) &= S v(t), \ v (0)= v_0, \; t\geq 0\\
\end{split}
\end{equation}
where $S \in R^{q \times q}$ is a constant matrix.

The regulated output for each subsystem is defined as
\begin{equation}\label{}
e_i(t) = y_i(t) - y_0(t),\ i=1,\dots,N,
\end{equation}
where $y_0(t)=-Fv(t)$.

Let $\mathcal{C}([-r,0],R^w)$ with  $r >0$ be
the Banach space of continuous functions mapping the interval $[-r,0]$
into $R^w$ endowed with the supremum norm. We assume $ x_{i0} \in
\mathcal{C}([-\tau_{con},0],R^n)$.

The plant \eqref{sys} and \eqref{exo} can be viewed as a multi-agent systems with the exosystem \eqref{exo} as the leader and the $N$ subsystems of \eqref{sys} as the followers. The communication
topology can be described by a directed graph $\mathcal{\bar{G}}=(\mathcal{\bar{V}},\mathcal{\bar{E}})$\footnote{See Appendix for a summary of digraph.}, where $\mathcal{\bar{V}}=\{0,1,\dots,N\}$ is the node set with the node 0 associated with the exosystem \eqref{exo} and all the other nodes associated with the $N$ subsystems \eqref{sys}, and $\mathcal{\bar{E}}$ is the edge set. The edge $(j,i) \in \mathcal{\bar{E}},\ i\neq j,\ i,j=0,\dots,N$, if and only if the control $u_i,\ i=1,\dots,N,$ can access the state $x_j$ and / or the output $y_j$ of subsystem $j,\ j=0,\dots,N$.
If $(j,i) \in \mathcal{\bar{E}}$, node $j$ is called a neighbor of the node $i$. We use $\mathcal{{N}}_i$ to denote the neighbor set of node $i$ with respect to $\mathcal{\bar{V}}$.

Due to the communication constraint and the communication time-delay,  we are limited to consider the class of distributed control laws with the communication delay.
Mathmetically, such a control law is described as follows:
\begin{equation}\label{ctr}
 \begin{split}
u_i(t)= & k(z_i(t),z_j(t),x_i(t-\tau_{com}),x_j(t-\tau_{com}),j\in \mathcal{{N}}_i)\\
\dot{z}_i(t)=& g(z_i(t),y_i(t-\tau_{com}),y_j(t-\tau_{com}),j\in \mathcal{{N}}_i)\\
z_i(\theta) =&   z_{i0}(\theta), \   \theta \in [-\tau_{com},0], i=1,\dots,N
 \end{split}
\end{equation}
where $z_i \in R^{n_{z}}$, $ z_{i0} \in \mathcal{C}([-\tau_{com},0],R^{n_{z}})$, $k$ and $g$ are linear functions of their arguments, $\tau_{com}\geq 0$ represents the communication delay among the agents. The control law \eqref{ctr} is called a distributed dynamic state feedback control law, and is further called a distributed dynamic output feedback control law if the function $k$ is independent of any state variable.

In recent years, the cooperative output regulation problem of multi-agent systems has received extensive attention \cite{Su1,Su5,Su3,Jiang1}. The problem is interesting because its formulation includes the leader-following consensus, synchronization or formation as special cases. Like the output regulation problem of a single linear system \cite{Da,Fr,FW}, there are two approaches to handling the cooperative output regulation problem of multi-agent systems. The first one is called feedforward design
\cite{Su1,Su5}. This approach makes use of the solution of the regulator
equations and a distributed observer to design an appropriate feedforward term to exactly cancel
the steady-state  tracking error. The second one is called
distributed internal model design \cite{Su3,Jiang1}. This approach employs
a distributed internal model
to convert the cooperative output regulation problem of an uncertain multi-agent system to a
simultaneous eigenvalue assignment problem of a multiple augmented system composed of the
given multi-agent system and the distributed internal model. The internal
model approach has at least two advantages over the feedforward design approach in that it can tolerate
perturbations of the plant parameters, and it does not need to solve the regulator equations.

More recently, the feedforward approach was further extended to the cooperative output regulation problem for exactly known linear multi-agent systems with time-delay \cite{cns}.
However, since this approach cannot handle the model uncertainties and the control law has to rely on the solution to the regulator equations, we will further develop a distributed internal
model approach to deal with the cooperative output regulation problem of uncertain multi-agent systems subject to both input delay and communication delay.

As a special case of the cooperative output regulation problem,  the leader-following consensus problem of linear multi-agent systems has been studied in several papers.
Some typical references that handle the communication time-delay are \cite{HH}, \cite{Jia}, \cite{Moreau}, \cite{Gao}, \cite{Saber}, \cite{Tian1}, \cite{Long} and \cite{Cheng}. In particular, in \cite{Saber},  the communication time-delays were considered in the leaderless consensus problem for single-integrator multi-agent systems
under undirected and fixed network topology.
In \cite{Cheng}, the leader-following consensus problem of double integrator multi-agent systems with non-uniform time-varying communication delays was studied under fixed and switching topologies.
On the other hand, input delay is also inevitable due to the processing and connecting time for packets arriving at each agent \cite{Zhoua}. Cooperative control of multi-agent systems with input delay has been studied in, say, \cite{Tian}, \cite{Xie}, \cite{Zhoua} and the references therein. In particular, the reference \cite{Zhoua} considered the leaderless consensus problem of high-order linear multi-agent systems with both communication delay and input delay with directed and fixed network.

As mentioned before, the problem formulation of this paper is general enough to include the leader-following consensus problem of  general multi-agent systems with
both communication delay and input delay as a special case. Moreover,  by adopting the distributed internal model approach, our control law is able to
handle model uncertainty, and simultaneously achieve asymptotic tracking and disturbance rejection for a large class of signals  generated by a linear autonomous system called exosystem.

Technically, this paper is most relevant to \cite{ijc} and \cite{Su3}. Specifically, reference  \cite{ijc} studied a special case of this paper with $N = 1$ in the system (\ref{sys}).
For this case, since there is no communication constraint on the control law (\ref{ctr}), we can use the full state feedback control or the full output feedback control to handle the problem. However, in the current case, we have to employ distributed control law which makes the design of our control law much more complicated. On the other hand, reference  \cite{Su3} treated
the same problem as this paper for a special case of the system (\ref{sys}) with
$\tau_{con} = 0$ by a special case of the control law (\ref{ctr}) with $\tau_{com} = 0$. However, due to the input delay and communication delay, the proof of the main results of this paper  is much more sophisticated than the proof of the main results in \cite{Su3}. We have to introduce or establish some specific technical lemmas to establish our main results.

The rest of this paper is organized as follows.
Section \ref{spfap} gives the problem formulation and some preliminaries. A general framework is established in Section \ref{sagf}.
Section \ref{smr} presents our main results. One example is used in Section \ref{se} to illustrate our results.
Finally, we close the paper with some concluding remarks in Section \ref{sc}.

{\bf Notation.} For $x_i\in R^{n_i}$, $i=1,...,m$, col$(x_1,\cdots,x_m)=[x_1^T,\cdots,x_m^T]^T$. For any matrix $X\in R^{n\times m}$, $\text{vec}(X)=\text{col}(X_1,\dots,X_m)$ where $X_i,i=1,\dots,m$, is the $i^{th}$ column of $X$. $\otimes$ denotes the Kronecker product of matrices. Let $\mathcal {C}$ denote the complex plane. For $\lambda \in \mathcal {C}$, let $Re\{\lambda \} >0$ denote the real part of $\lambda$.

\section{Problem formulation and preliminaries}\label{spfap}

Like in \cite{Su3}, all matrices in \eqref{sys} can be uncertain. Let $\bar{A}_i=A_i+\delta_i A,\bar{B}_i=B+\delta_i B,
\bar{E}_i=E_i+\delta_i E,\bar{C}_i=C_i+\delta_i C$, where $A_i,B_i,E_i,C_i$ represent the
nominal part of these matrices, and $\delta_i A,\delta_i B,\delta_i E,\delta_i C$ are the perturbations of these matrices. For convenience, we denote the system uncertainties with a vector
\begin{equation*}
w = \left[
                                                                                                         \begin{array}{c}
                                                                                                           \text{vec}(\delta_1 A,\dots,\delta_N A)\\
                                                                                                           \text{vec}(\delta_1 B,\dots,\delta_N B)\\
                                                                                                           \text{vec}(\delta_1 E,\dots,\delta_N E)\\
                                                                                                           \text{vec}(\delta_1 C,\dots,\delta_N C)\\
                                                                                                         \end{array}
                                                                                                       \right] \in R^{Nn (n+m+p+q)}.
\end{equation*}

Now, we can state our problem as follows:
\begin{Definition}\label{Def1.1}
Linear cooperative robust output regulation problem: given the system \eqref{sys}, the exosystem \eqref{exo}, and a digraph $\mathcal{\bar{G}}$, design a control law of the form \eqref{ctr} such that the closed-loop system satisfies the properties \ref{Per1}
and \ref{Per2} as follows.

\begin{Property}\label{Per1}
The nominal closed-loop system is exponentially stable when $v=0$.
\end{Property}

\begin{Property}\label{Per2}
There exists an open neighborhood $W$ of $w=0$ such that, for any $w \in W$ and any initial conditions $x_{i0}$, $z_{i0}$ and $v_0$, the regulated output $\lim_{t \rightarrow \infty}e_i(t)= 0,\ i=1,\dots,N$.
\end{Property}

\end{Definition}

\begin{Remark}
It is noted that the problem studied in \cite{Su3} is a special case of the above problem when both
the communication delay $\tau_{com}$ and the input delay $\tau_{con}$ are zero. The presence of these two delays makes our problem formulation more realistic and, as will be seen later, the handling of the problem more challenging.
\end{Remark}

For the solvability of the above problem, some assumptions are stated as follows.

\begin{Assumption}\label{Ass1.0}
There exist matrices $A$, $B$, $C$ such that $A_i=A$, $B_i=B$, $C_i=C$, $i=1,\dots,N$.
\end{Assumption}

\begin{Assumption}\label{Ass1.1}
All the eigenvalues of $S$ are on the imaginary axis.
\end{Assumption}

\begin{Assumption}\label{Ass1.2}
The matrix pair $(A,B)$ is stabilizable.
\end{Assumption}

\begin{Assumption}\label{Ass1.3}
The matrix pair $(C,A)$ is detectable.
\end{Assumption}

\begin{Assumption}\label{Ass1.4}
For all $\lambda \in \sigma(S)$, where $\sigma(S)$ denotes the spectrum of $S$,
\begin{equation}\label{tzero}
\text{rank} \left(
                                      \begin{array}{cc}
                                          A -\lambda I_n &  B \\
                                          C & 0 \\
                                      \end{array}
                                    \right)=n+p.
\end{equation}
\end{Assumption}

\begin{Assumption}\label{Ass1.5}
The digraph $\bar{\mathcal{G}}$ contains a directed spanning tree with the node $0$ as the root.
\end{Assumption}

\begin{Assumption}\label{Ass1.21}
$A$ has no eigenvalues with positive real parts.
\end{Assumption}

\begin{Remark}
Assumptions \ref{Ass1.0} to \ref{Ass1.5} are standard ones and they are needed in \cite{Su3} even if there are no communication delay and input delay.
And Assumption \ref{Ass1.21} is additional and it is made so that the delayed system can be stabilized by using the low gain method introduced in \cite{Zhou3}.
\end{Remark}

\section{A general framework}\label{sagf}

To construct a specific control law, let $\mathcal{\bar{A}}=\left[ a_{ij}\right] \in R^{(N+1)\times(N+1)}$ and $\mathcal{\bar{L}}=[l_{ij}] \in R^{(N+1)\times(N+1)}$ be the
 weighted adjacent matrix and Laplacian of the digraph $\mathcal{\bar{G}}$,  respectively.   Let $\Delta$ be an $N \times N$ nonnegative diagonal matrix whose $i^{th}$ diagonal element is $a_{i0}$. Then, we have \cite{HH,Su1} \begin{equation*}
\mathcal{\bar{L}} = \left(
       \begin{array}{c|c}
         0   & 0_{1\times N} \\ \hline
         -\Delta  \mathbf{1}_N  & H \\
       \end{array}
     \right)
\end{equation*}
where $ \mathbf{1}_N $ is an $N \times 1$ column vector whose elements are all $1$ and $H$ satisfies $H \mathbf{1}_N =\Delta  \mathbf{1}_N$.

In terms of the elements of $\mathcal{\bar{A}}$, we can define a virtual regulated output $e_{vi}(t)$ for each follower subsystem $i$ as follows:
\begin{equation}\label{sysvirt}
e_{vi}(t) =\sum_{j\in \mathcal{{N}}_i} a_{ij} (y_i(t)-y_j(t)),\ i=1,\dots,N.
\end{equation}
Note that the subsystem $e_{vi}(t)$ can access the regulated error $(y_i(t)-y_j(t))$ if and only if the node $j$ is the neighbor of the node $i$.

\begin{Remark}\label{Remeig}
Let $e=\mbox{col}(e_1,\dots,e_{N})$ and $e_v=\mbox{col}(e_{v1},\dots,e_{vN})$. Then it can be verified that $e_v = (H \otimes I_p) e$.  By Lemma 4 of  \cite{HH} or  Lemma 1 of \cite{Su1},
the matrix $-H$ is Hurwitz if and only if Assumption \ref{Ass1.5} is satisfied. Thus, under Assumption \ref{Ass1.5}, $e_v = 0$ iff $e = 0$.
\end{Remark}

In order to  make use of the internal model principle  to handle the systems with input delay and communication delay, we need to generalize  the concept of the minimum p-copy internal model
to the following form:

\begin{Definition}\label{Definter}
A pair of matrices $(G_1,G_2)$ is said to be the minimal p-copy
internal model of the matrix $S$ if the pair takes the following
form:
\begin{equation}\label{inter}
 {G}_1=\text{block diag} \underbrace{(\beta,\dots,\beta)}_{p-tuple},\ {G}_2=\text{block diag} \underbrace{(\sigma,\dots,\sigma)}_{p-tuple}
\end{equation}
where $\beta$ is a constant square matrix whose characteristic polynomial equals the minimal polynomial of $S$, and $\sigma$ is a constant column vector such that $(\beta,\sigma)$ is controllable.
\end{Definition}

Having defined the virtual regulated output $e_{vi}(t)$ and introduced the p-copy internal model, we can describe our distributed dynamic state feedback control law
as follows:
\begin{equation}\label{ctr1}
\begin{split}
u_i(t)=&  K_x \eta_i(t) + K_z z_i(t)\\
\dot{z}_i(t) =&  {G}_1 z_i(t) + {G}_2 e_{vi}(t-\tau_{com}),\ i=1,\dots,N  \\
\end{split}
\end{equation}
where $\eta_i(t)= \sum_{j\in\mathcal{{N}}_i} a_{ij}\left(x_i\left(t-\tau_{com}\right)-x_j\left(t-\tau_{com}\right)\right) $, $x_{0}(t)=0$, $z_i(t) \in R^{n_z}$ with $n_z$ to be specified later, $(K_x,K_z)$ are constant matrices of appropriate dimensions to be designed later, $(G_1,G_2)$ are defined in \eqref{inter}, and, respectively,
our distributed dynamic output feedback control law as follows:
\begin{equation}\label{ctr2}
\begin{split}
u_i(t)=&   K_1 z_i(t) + K_{2} \hat{\eta}_i(t) \\
\dot{z}_i(t) =&  {G}_1 z_i(t) + {G}_2 e_{vi}(t-\tau_{com})  \\
\dot{\xi}_i(t) =& A \xi_i(t) + Bu_i(t-\tau)  -LC \hat{\eta}_i(t) + Le_{vi}(t-\tau_{com}), \ i=1,\dots,N
\end{split}
\end{equation}
where $\xi_i(t) \in R^n$, $\hat{\eta}_i(t)= \sum_{j\in\mathcal{{N}}_i}a_{ij}\left(\xi_i(t)-\xi_j(t)\right) $, $\xi_0(t)=0$, and $z_i(t) \in R^{n_z}$ with $n_z$ to be specified later, $\tau=\tau_{com}+\tau_{con}$, $(K_1,K_{2},L)$ are constant matrices of appropriate dimensions to be designed later and $(G_1,G_2)$ are defined in \eqref{inter}.

Let $x=\mbox{col}(x_1,\dots,x_N)$,  $z= \mbox{col}(z_1,\dots,z_N)$, $\xi= \mbox{col}(\xi_1,\dots,\xi_N)$, $u=\mbox{col}(u_1,\dots,u_N)$, $\tilde{A}=\text{blockdiag}(\bar{A}_1,$ $\dots,\bar{A}_N)$, $\tilde{B}=\text{blockdiag}(\bar{B}_1,\dots,\bar{B}_N)$, $\tilde{E}=(\bar{E}_1^T,\dots,$ $\bar{E}_N^T)^T$, $\tilde{C}=(H\otimes I_p)\text{block diag}(\bar{C}_1,\dots,\bar{C}_N)$, $\tilde{F}=(\Delta \textbf{1}_N)\otimes F$, $\tilde{L}=I_N\otimes L$, $\tilde{G}_1=I_N\otimes G_1$, $\tilde{G}_2=I_N\otimes G_2$.
Then, we define an auxiliary system as follows:
\begin{equation}\label{sysauxi}
\begin{split}
\dot{x}(t) &= \tilde{A} x(t) +  \tilde{B} u(t-\tau_{con}) + \tilde{E} v(t), t \geq 0\\
\dot{v}(t) &= S v(t), t \geq 0\\
e_v(t) &= \tilde{C} x(t) + \tilde{F} v(t), t \geq 0.\\
\end{split}
\end{equation}

Clearly, the matrix pair $(\tilde{{G}}_1, \tilde{{G}}_2)$ is the minimal pN-copy internal model of the matrix $S$. Thus, by Definition \ref{Definter}, the following system
\EQ \label{imauxi}
\dot{z}(t)  =  \tilde{ {G}}_1 z (t) +  \tilde{ {G}}_2 e_v(t-\tau_{com}) , t \geq 0
\EN
 is an internal model of (\ref{sysauxi}).
The composition of the auxiliary system (\ref{sysauxi}) and the (\ref{imauxi}) is called the augmented system of (\ref{sysauxi}) and is put as follows:
\begin{equation}\label{syscomp}
\begin{split}
\dot{x}(t) &= \tilde{A} x(t) +  \tilde{B} u(t-\tau_{con}) + \tilde{E} v(t), t \geq 0\\
\dot{z}(t) & =  \tilde{ {G}}_1 z (t) +  \tilde{ {G}}_2 e_v(t-\tau_{com}), t \geq 0 \\
\dot{v}(t) &= S v(t), t \geq 0\\
e_v(t) &= \tilde{C} x(t) + \tilde{F} v(t), t \geq 0.
\end{split}
\end{equation}

\begin{Remark}\label{Lemi}
It can be seen that the internal model in \cite{ijc} is a special case of \eqref{imauxi} by setting $\tau_{com} = 0$. It is shown in  Lemma 1.27 of \cite{J.H} that if the matrix pair $(G_1,G_2)$ is the minimal p-copy internal
model of the matrix $S$, then the following matrix equation
\begin{equation}\label{ire3}
X S - G_1 X = G_2 Y
\end{equation}
has a solution $X $ only if $Y=0$. This property is the key for establishing the following result.
\end{Remark}

The role of an internal model is to convert the output regulation problem of the given plant \eqref{sysauxi} to the stabilization problem
of the augmented system (\ref{syscomp}). To be more precise, we have the following lemma.
\begin{Lemma}\label{Remsta}
Under Assumption \ref{Ass1.1}, \\
(i) suppose a static state feedback control law of the form
\begin{equation}\label{Ctr1n01}
u(t)=\tilde{K}_x x(t-\tau_{com})+ \tilde{K}_{z} {z}(t)
\end{equation}
stabilizes the nominal plant of the augmented system \eqref{syscomp}. Then, the dynamic state feedback control law of the form
\begin{equation}\label{Ctr1n1}
\begin{split}
u(t)=& \tilde{K}_x x(t-\tau_{com})+ \tilde{K}_{z} z (t)\\
\dot{z}(t) =& {\tilde G}_1 z (t) +{\tilde G}_2 e_v(t-\tau_{com})  \\
\end{split}
\end{equation}
solves the robust output regulation problem of the auxiliary system (\ref{sysauxi}).\\
(ii) suppose a dynamic output feedback control law of the form
\begin{equation}\label{Ctr1n02}
\begin{split}
u(t)=& \tilde{K}_1 z(t) + \tilde{K}_{2} \xi(t)\\
\dot{\xi}(t)=& \tilde{S}_{1}\xi(t)+ \tilde{S}_2 e_v(t-\tau_{com}) + \tilde{S}_{3} \zeta(t- \tau )
\end{split}
\end{equation}
where $ \zeta=\mbox{col}(z,\xi)$,
stabilizes the nominal plant of the augmented system \eqref{syscomp}. Then, the dynamic output feedback control law of the form
\begin{equation}\label{Ctr1n2}
\begin{split}
u(t)=& \tilde{K}_1 z(t) + \tilde{K}_{2} \xi(t)\\
\dot{z}(t) =& \tilde{ {G}}_1 z (t) +  \tilde{ {G}}_2 e_v(t-\tau_{com}) \\
\dot{\xi}(t)=& \tilde{S}_{1}\xi(t)+ \tilde{S}_2 e_v(t-\tau_{com}) + \tilde{S}_{3} \zeta(t- \tau )
\end{split}
\end{equation}
solves the robust output regulation problem of the auxiliary system (\ref{sysauxi}).
By Remark \ref{Remeig}, under Assumption \ref{Ass1.5}, either of the two control laws also solves the cooperative robust output regulation problem of the given plant (\ref{sys}).
\end{Lemma}

Before giving the proof of Lemma \ref{Remsta}, we still need some remarks. First, under the coordinate transformation $z(t)=\bar{z}(t-\tau_{com}), \xi(t)=\bar{\xi}(t-\tau_{com})$,  the closed-loop system composed of system \eqref{sysauxi} and \eqref{Ctr1n1} or \eqref{Ctr1n2} can be put into the following form:
\begin{equation}\label{cl3}
\begin{split}
\dot{x}_c(t) &= \sum_{i=0}^{1} A_{cwi} x_c(t-\bar{\tau}_i) + B_{cw} v(t)\\
e_v(t) &=  \sum_{i=0}^{1} C_{cwi} x_c(t-\bar{\tau}_i) + D_{cw} v(t) \\
\dot{v}(t) &=Sv(t) \\
\end{split}
\end{equation}
where $\bar{\tau}_0=0$, $\bar{\tau}_1=\tau$, under the dynamic state feedback, $x_c=\mbox{col}(x,\bar{z})$, and
\begin{equation*}\label{}
\begin{split}
A_{cw0} =& \left(
          \begin{array}{cc}
            \tilde{A}  & 0 \\
            \tilde{G}_2 \tilde{C}  & \tilde{G}_1 \\
          \end{array}
        \right), \
A_{cw1} =  \left(
          \begin{array}{cc}
            \tilde{B} \tilde{K}_x & \tilde{B} \tilde{K}_z \\
            0 & 0  \\
          \end{array}
        \right) \\
B_{cw}= & \left(
            \begin{array}{c}
              \tilde{E} \\
            \tilde{G}_2 \tilde{F}  \\
            \end{array}
          \right) , \
C_{cw0} =  \left(
             \begin{array}{cc}
               \tilde{C} & 0\\
             \end{array}
           \right) \\
C_{cw1} =& \left(
             \begin{array}{cc}
               0 & 0\\
             \end{array}
           \right) ,  \
D_{cw}=  \tilde{F}  \\
\end{split}
\end{equation*}
and, under the dynamic output feedback, $x_c = \mbox{col}(x ,\bar{z},\bar{\xi} )$, and
\begin{equation*}\label{}
\begin{split}
A_{cw0} =& \left(
          \begin{array}{ccc}
            \tilde{A} & 0 & 0\\
            \tilde{G}_2 \tilde{C}  & \tilde{G}_1 & 0\\
            \tilde{S}_2 \tilde{C} & 0 & \tilde{S}_1 \\
          \end{array}
        \right) ,\
A_{cw1} = \left(
          \begin{array}{cc}
             0 &
                    \begin{array}{cc}
                      \tilde{B} \tilde{K}_1 & \tilde{B} \tilde{K}_2 \\
                    \end{array}
                    \\
                \begin{array}{c}
                  0 \\
                  0 \\
                \end{array}    &
                                \begin{array}{c}
                                  0 \\
                                  \tilde{S}_3 \\
                                \end{array}\\
          \end{array}
        \right) \\
B_{cw}= & \left(
            \begin{array}{c}
              \tilde{E} \\
              \tilde{G}_2 \tilde{F}  \\
              \tilde{S}_2 \tilde{F}  \\
            \end{array}
          \right)  , \
C_{cw0} = \left(
             \begin{array}{ccc}
               \tilde{C}  & 0 & 0\\
             \end{array}
           \right) \\
C_{cw1} =& \left(
             \begin{array}{ccc}
               0 & 0 & 0\\
             \end{array}
           \right) , \
D_{cw}=   \tilde{F} \\
\end{split}
\end{equation*}

\begin{Remark}\label{lemCyber} It can be deduced from Lemma 2.1 of \cite{ijc}, under Assumption \ref{Ass1.1}, if  the closed-loop system
\eqref{cl3} satisfies Property \ref{Per1}, then, for each $w\in W$,
and any matrix $B_{cw}$ of appropriate dimension, there exists a
unique matrix $X_{cw}$ that satisfies the following matrix equation:
\begin{equation}\label{lem51}
X_{cw} S = \sum_{i=0}^{1} A_{cwi} X_{cw} e^{-S \bar{\tau}_i} + B_{cw}.
\end{equation}
Moreover, by Lemma 2.2 of \cite{ijc}, under Assumption \ref{Ass1.1}, if the
controller \eqref{Ctr1n1} or \eqref{Ctr1n2} renders the closed-loop system
\eqref{cl3} Property \ref{Per1}, then, the same controller
solves the linear robust output regulation problem if and only if,
for each $w\in W$, there exists a unique matrix $X_{cw}$ that
satisfies the following matrix equations:
\begin{equation}\label{re1}
\begin{aligned}
X_{cw} S &=  \sum_{i=0}^{1} A_{cwi} X_{cw} e^{-S \bar{\tau}_i} + B_{cw}   \\
0   &=  \sum_{i=0}^{1} C_{cwi} X_{cw} e^{-S\bar{\tau}_i} +  D_{cw} .\\
\end{aligned}
\end{equation}
\end{Remark}

Now, we will give the proof of Lemma \ref{Remsta} as follows.

\begin{Proof}
Note that the closed-loop system \eqref{cl3} can also be viewed as a composition
of the augmented system \eqref{syscomp} and a static state feedback control of the form $u(t)=\tilde{K}_x x(t-\tau_{com})+ \tilde{K}_{z} {z}(t)$
 $($respectively, a dynamic output feedback control law of the form $u(t)= \tilde{K}_1 z(t) + \tilde{K}_{2} \xi(t), \dot{\xi}(t)= \tilde{S}_{1}\xi(t)+ \tilde{S}_2 e_v(t-\tau_{com}) + \tilde{S}_{3} \zeta(t- \tau ), \text{where} \ \zeta=\mbox{col}(z,\xi)$$)$. Thus, the closed-loop system \eqref{cl3} satisfies Property \ref{Per1}. By Remark \ref{lemCyber}, under Assumption \ref{Ass1.1},
 it suffices to prove that the matrix equations \eqref{re1} have a unique solution $X_{cw}$ under either the static state feedback controller
 or the dynamic output feedback controller.
 In fact, by Remark \ref{lemCyber}, the first equation of \eqref{re1} has one unique solution $X_{cw}$. Thus, we only need to prove that $X_{cw}$ also satisfies the second equation of \eqref{re1}. We will do so for  the static state feedback control case and the dynamic output feedback case, respectively.

Part (i): Let
$X_{cw}=\left(
                                                                      \begin{array}{c}
                                                                        X_w \\
                                                                         Z_w \\
                                                                      \end{array}
                                                                    \right)$
with $X_{w} \in R^{Nn \times q} $ and expand the first  equation of
(\ref{re1}) to the following form:
\begin{equation}\label{ire1}
\begin{split}
X_w S=& \tilde{A} X_w +  \tilde{B} ( \tilde{K}_x X_w   + \tilde{K}_z Z_w )e^{-S  {\tau} } + \tilde{E} \\
Z_w S=& \tilde{G}_1 Z_w  + \tilde{G}_2 Y_w
\end{split}
\end{equation}
where
\begin{equation}\label{ire2}
Y_w = \tilde{C} X_w + \tilde{F}.
\end{equation} Since
the second equation of  (\ref{ire1}) is in the form (\ref{ire3}), by
Remark \ref{Lemi}, $Y_w = 0$. That is, $X_{cw}$ also satisfies the
second equation of (\ref{re1}).

Part (ii):
Let $X_{cw}=\left(
                                                                      \begin{array}{c}
                                                                        X_w \\
                                                                         {Z}_w \\
                                                                      \end{array}
                                                                    \right)$
with $X_w \in R^{Nn \times q}$, and ${Z}_w \in R^{N(n_{z}+n) \times q}$. Partition ${Z}_w$ to
$Z_{w}=\left(
                                                                      \begin{array}{c}
                                                                        \bar{Z}_w \\
                                                                        \hat{Z}_w \\
                                                                      \end{array}
                                                                    \right)$,
where $\bar{Z}_w \in R^{Nn_z \times q}$ with $Nn_z$ the dimension of $\tilde{G}_1$.
Then, it can be verified that, under the control law $u(t)= \tilde{K}_1 z(t) + \tilde{K}_{2} \xi(t), \dot{\xi}(t)= \tilde{S}_{1}\xi(t)+ \tilde{S}_2 e_v(t-\tau_{com}) + \tilde{S}_{3} \zeta(t- \tau ), \text{where} \ \zeta=\mbox{col}(z,\xi)$,
the first equation of (\ref{re1}) can be expanded to the following
form:
\begin{equation}\label{ire1o}
\begin{split}
X_w S=& \tilde{A} X_w + \tilde{B} ( \tilde{K}_1  \bar{Z}_w + \tilde{K}_2 \hat{Z}_w) e^{-S  {\tau} } + \tilde{E} \\
\bar{Z}_w  S=& \tilde{G}_1 \bar{Z}_w   + \tilde{G}_2 Y_w\\
\hat{Z}_w  S=&  \tilde{S}_{1}  \hat{Z}_w + \tilde{S}_2 Y_w +  \tilde{S}_{3} {Z}_w  e^{-S {\tau} }\\
\end{split}
\end{equation}
where
\begin{equation}\label{ire2o}
Y_w = \tilde{C} X_w + \tilde{F}.
\end{equation}
Since the second equation of  (\ref{ire1o}) is in the form (\ref{ire3}), by
Remark \ref{Lemi}, $Y_w = 0$. That is, $X_{cw}$ also satisfies the
second equation of (\ref{re1}).
\end{Proof}

\begin{Remark}\label{Remimp0}
In order to apply Lemma \ref{Remsta} to our problem,  it is not enough to show that the nominal part of the augmented system \eqref{syscomp} is stabilizable by a static state feedback control law of
the form (\ref{Ctr1n01}) or a dynamic output feedback control law of the form (\ref{Ctr1n02}). We actually need to show that the nominal part of the augmented system \eqref{syscomp} is stabilizable by a distributed static state feedback control law of the form  $u_i(t)=  K_x \eta_i(t) + K_z z_i(t)$, $i = 1, \dots, N$ (or a distributed dynamic output feedback control law of the form $u_i(t)=  K_1 z_i(t) + K_{2} \hat{\eta}_i(t),\dot{\xi}_i(t)= {S}_{1}\xi_i(t) +  {S}_2 e_{vi}(t-\tau_{com}) + {S}_{3} \zeta_i(t- \tau )+ S_{4} \hat{\eta}_i (t) + S_{5} \hat{\eta}_i (t-\tau) $,where $\zeta_i=\mbox{col}(z_i,\xi_i)$, $i = 1, \dots, N$). As a result, the distributed state feedback control law
(\ref{ctr1}) ( or the distributed output feedback control law
(\ref{ctr2})) solves the cooperative output regulation problem of the system (\ref{sys}).
 What makes this stabilization problem much more challenging than the problem in \cite{Su3} is that the augmented system \eqref{syscomp} is subject to both input delay and communication delay. We need to first establish a few lemmas to lay the foundation of our approach.

\end{Remark}

%

\section{Main result}\label{smr}

To establish some Lemmas in this section, we need to first cite the following lemma.
\begin{Lemma}(Lemma 3.2 in \cite{cns})\label{Lem2}
Consider the system
\begin{equation}\label{sw1}
\begin{split}
\dot{\zeta}(t) &= M_0 \zeta(t) + \sum_{i=1}^{p} M_i \zeta(t-\Delta_i) + N \xi(t)\\
\zeta(\theta)&= \zeta_0(\theta), \ \theta \in [-\Delta,0]
\end{split}
\end{equation}
where $M_i \in R^{n\times n},\ i=0,1,\dots,p$, $N \in R^{n\times m}$
are some constant matrices,
$0<\Delta_1<\Delta_2<\dots<\Delta_p=\Delta$ are arbitrary time-delays, $\zeta_0 \in
\mathcal{C}([-\Delta,0],R^n)$, and $\xi(t)$ is any measurable,
essentially bounded function over $[0, \infty)$. Assume that the
origin of the unforced $\zeta(t)$ system is exponentially stable
and $\lim_{t\rightarrow \infty}\xi(t)=0$. Then, $\lim_{t\rightarrow
\infty}\zeta(t)=0$. Moreover, $\lim_{t\rightarrow \infty}\zeta(t)=0$
exponentially if $\lim_{t\rightarrow \infty}\xi(t)=0$ exponentially.
\end{Lemma}

\begin{Lemma}\label{LemAB}
Suppose that Assumptions \ref{Ass1.1}, \ref{Ass1.2}, \ref{Ass1.4} and \ref{Ass1.21} are satisfied.
Consider the system of the form
\begin{equation}\label{syslem1}
\dot{x}_i(t) = \mathcal{A} x_i(t) + \lambda_i \mathcal{B} u_i(t-\tau),\ i=1,\dots,N
\end{equation}
where $x_i \in R^{(n+n_z)}$, $u_i \in R^{m}$, $\mathcal{A}=\left(
       \begin{array}{cc}
         A & 0 \\
         G_2C & G_1 \\
       \end{array}
     \right)
$, $\mathcal{B}=\left(
                 \begin{array}{c}
                   B \\
                   0 \\
                 \end{array}
               \right)
$, and $\lambda_i \in \mathcal {C}$ with $Re\{\lambda_i\} >0$. Then, there exists a matrix $K \in R^{m \times (n+n_z)}$ such that the state feedback control law $u_i(t)=Kx_i(t)$,
$i=1,\dots,N$, asymptotically stabilize all subsystems of the system \eqref{syslem1}.

\end{Lemma}

\begin{Proof}
Under Assumptions \ref{Ass1.1}, \ref{Ass1.2} and \ref{Ass1.4}, by Lemma 1.26 of \cite{J.H}, $(\mathcal{A},\mathcal{B})$ is stabilizable. Moreover, under additional Assumption \ref{Ass1.21}, we have that $\mathcal{A}$ has no eigenvalues with positive real parts. Therefore, there exists a nonsingular matrix $ \mathcal {T}$ such that
\begin{equation}\label{}
\bar{\mathcal{A}}=\mathcal {T} \mathcal{A} \mathcal {T}^{-1}=\left(
          \begin{array}{cc}
            \mathcal{A}_1 & 0 \\
            0 & \mathcal{A}_2 \\
          \end{array}
        \right),\
\bar{\mathcal{B}}= \mathcal {T}\mathcal{B}=\left(
              \begin{array}{c}
                \mathcal{B}_1 \\
                \mathcal{B}_2 \\
              \end{array}
            \right)
\end{equation}
where all the eigenvalues of the matrix $\mathcal{A}_2$ have negative real parts, all the eigenvalues of the matrix $\mathcal{A}_1$ are on the imaginary axis and $(\mathcal{A}_1,\mathcal{B}_1)$ is controllable. Then, system \eqref{syslem1} is equivalent to the following system:
\begin{equation}\label{eqAB1}
\begin{split}
\dot{\chi}_{i1}(t) =& \mathcal{A}_1 \chi_{i1}(t) + \lambda_i \mathcal{B}_1 u_i(t-\tau) \\
\dot{\chi}_{i2}(t) =& \mathcal{A}_2 \chi_{i2}(t) + \lambda_i \mathcal{B}_2 u_i(t-\tau),\ i=1,\dots,N.\\
\end{split}
\end{equation}

By Lemma 1 of \cite{Zhoua}, there exists a matrix $\bar{K}_1=- \nu_1^{-1} \mathcal{B}_1^T P e^{\mathcal{A}_1\tau}$, where $\nu_1 \in R$ satisfies
\begin{equation}\label{}
0< \nu_1 \leq Re(\lambda_i),\ i=1,\dots,N,
\end{equation}
and $P$ is the positive definite solution of the ARE
\begin{equation}\label{AR}
\mathcal{A}_1^T P + P \mathcal{A}_1 - P\mathcal{B}_1\mathcal{B}_1^TP=-\gamma P
\end{equation}
with some sufficiently small $\gamma >0$ such that,  for $i=1,\dots,N$, the systems $\dot{\chi}_{i1}(t) = \mathcal{A}_1 \chi_{i1}(t) + \lambda_i \mathcal{B}_1 \bar{K}_1 \chi_{i1}(t-\tau)$ are all asymptotically stable.

Let $\bar{K} =(\bar{K}_1,0)$. Then, under the control law $u_i(t)=\bar{K} \chi_i(t)$, the closed-loop system of \eqref{eqAB1} is as follows.
\begin{equation}\label{eqAB2}
\begin{split}
\dot{\chi}_{i1}(t) =& \mathcal{A}_1 \chi_{i1}(t) + \lambda_i \mathcal{B}_1 \bar{K}_1 \chi_{i1}(t-\tau) \\
\dot{\chi}_{i2}(t) =& \mathcal{A}_2 \chi_{i2}(t) + \lambda_i \mathcal{B}_2 \bar{K}_1 \chi_{i1}(t-\tau),\ i=1,\dots,N\\
\end{split}
\end{equation}

Since for $i=1,\dots,N$, $\chi_{i1}(t)$ subsystem is asymptotically stable, by Lemma \ref{Lem2}, for $i=1,\dots,N,$ $\chi_{i2}(t)$ subsystem is asymptotically stable. The proof is thus completed with $K=\bar{K}\mathcal{T}$.
\end{Proof}

\begin{Lemma}\label{Lem1sta}
Consider the system of the form
\begin{equation}\label{linm1}
\begin{aligned}
\dot{x}_c(t) =& \left(
               \begin{array}{cc}
                  I_N\otimes A  & 0 \\
                  H\otimes G_2C  &  I_N\otimes G_1  \\
               \end{array}
             \right) x_c(t) +\left(
                              \begin{array}{cc}
                                H\otimes B & I_N\otimes B \\
                                0_{N n_z \times Nm} & 0_{N n_z \times Nm} \\
                              \end{array}
                            \right)u_c(t-\tau)\\
x_c (\theta) =&  x_{c0}(\theta), \ \theta \in [-\tau,0]
\end{aligned}
\end{equation}
where $x_c \in R^{N(n+n_z)}$, $u_c \in R^{2Nm}$,  $(G_1,G_2)$ is the minimal p-copy internal model of $S$ as defined in \eqref{inter},
and $x_{c0} \in \mathcal{C}\big([-\tau,0]$ $,R^{N(n+n_z)}\big)$. Then, under Assumptions \ref{Ass1.1}, \ref{Ass1.2}, \ref{Ass1.4} and \ref{Ass1.21}, there exist matrices $K_x \in R^{m \times n}$ and $K_z \in R^{m \times n_z}$, such that under the state feedback control law
$u_c(t)=\left(
          \begin{array}{cc}
            I_N \otimes K_x & 0_{Nm \times N n_z} \\
            0_{Nm \times N n} & I_N \otimes K_z \\
          \end{array}
        \right)x_c(t)$, system \eqref{linm1} is asymptotically stable if and only if Assumption \ref{Ass1.5} is satisfied.
\end{Lemma}

\begin{Proof}
(If Part:) Denote the eigenvalues of $H$ by $\lambda_i,\ i=1,\dots,N$. Under Assumption \ref{Ass1.5}, by Remark \ref{Remeig}, for $i=1,\dots,N$,   $\lambda_i$ have positive real parts.
Let $T_1$ be the nonsingular matrix such that $J_H=T_1HT_1^{-1}$ is in the Jordan form of $H$. Let
$T_2=\left(
\begin{array}{cc}
T_1\otimes I_n & 0 \\
0 & T_1\otimes I_{n_z} \\
\end{array}
\right)
$ and $\bar{x}_c(t)=T_2x_c(t)$. Then, $\bar{x}_c(t)$ is governed by the following system:
\begin{equation}\label{linm2}
\begin{aligned}
\dot{\bar{x}}_c(t) =& \left(
               \begin{array}{cc}
                 I_N\otimes A & 0 \\
                 J_H\otimes G_2C & I_N\otimes G_1 \\
               \end{array}
             \right) \bar{x}_c(t)  +
             \left(
               \begin{array}{cc}
                 J_H\otimes B & I_N\otimes B \\
                 0_{N n_z \times Nm} & 0_{N n_z \times Nm} \\
               \end{array}
             \right) \bar{T}_2 u_c (t-\tau)\\
\end{aligned}
\end{equation}
where $\bar{T}_2=\left(
                   \begin{array}{cc}
                     T_1\otimes I_m & 0 \\
                     0 & T_1\otimes I_m \\
                   \end{array}
                 \right)
$.

Denote ${u}_c=\mbox{col}( {u}_{c1},\dots, {u}_{cN})$ with $u_{ci} \in R^{2m}$.
Partition $I_{N(n+n_z)}$ as $I_{N(n+n_z)}=(M_1^T,\dots,M_N^T,Q_1^T$ $,\dots,Q_N^T)^T$, where $M_i \in R^{n \times (N(n+n_z))}$ and $Q_i \in R^{n_z \times (N(n+n_z))}$ and $I_{2Nm}$ as  $I_{2Nm}=(\bar{M}_1^T,\dots,\bar{M}_N^T,\bar{Q}_1^T$ $,\dots,\bar{Q}_N^T)^T$, where $\bar{M}_i \in R^{m \times (2Nm)}$ and $\bar{Q}_i \in R^{m \times (2Nm)}$.

Let $T_3=(M_1^T,Q_1^T,M_2^T,Q_{2}^T,\dots,M_N^T,Q_{N}^T)^T$, $\hat{x}_c=T_3\bar{x}_c$ and $\hat{x}_c=\mbox{col}(\hat{x}_{c1},\dots,\hat{x}_{cN})$. Then, the system \eqref{linm2} becomes
a lower triangular system whose diagonal blocks are of the form
\begin{equation}\label{linm3}
\begin{aligned}
\dot{\hat{x}}_{ci}(t) =& \left(
               \begin{array}{cc}
                  A & 0 \\
                 \lambda_i G_2C & G_1 \\
               \end{array}
             \right) \hat{x}_{ci}(t)   +
             \left(
               \begin{array}{cc}
                 \lambda_i B & B \\
                 0_{n_z \times m} & 0_{n_z \times m} \\
               \end{array}
             \right)   \hat{u}_{ci} (t-\tau) , \ i=1,\dots,N\\
\end{aligned}
\end{equation}
where $\hat{x}_{ci} \in R^{(n+n_z)}$, $\hat{u}_{ci} \in R^{2m}$, $\hat{u}_{c}=\mbox{col}(\hat{u}_{c1},\dots,\hat{u}_{cN})$, $\hat{u}_{c}=\bar{T}_3 \bar{T}_2 u_{c}$, and $\bar{T}_3=(\bar{M}_1^T,\bar{Q}_1^T,\bar{M}_2^T,\bar{Q}_{2}^T,\dots$ $,\bar{M}_N^T,\bar{Q}_{N}^T)^T$.

Let
$T_{4i}=\left(
              \begin{array}{cc}
                I_n & 0 \\
                0 & \lambda_i^{-1}I_{n_z} \\
              \end{array}
            \right)
$, and $\tilde{x}_{ci}(t)=T_{4i}\hat{x}_{ci}(t)$. Then, we get, for $i=1,\dots,N$,
\begin{equation}\label{linm4}
\begin{aligned}
\dot{\tilde{x}}_{ci}(t) =& \left(
               \begin{array}{cc}
                  A & 0 \\
                  G_2C & G_1 \\
               \end{array}
             \right) \tilde{x}_{ci}(t) + \lambda_i
             \left(
               \begin{array}{cc}
                 B & B \\
                 0 & 0 \\
               \end{array}
             \right) \bar{T}_{4i} \hat{u}_{ci} (t-\tau) \\
\end{aligned}
\end{equation}
where $\bar{T}_{4i}=\left(
                      \begin{array}{cc}
                        I_m & 0 \\
                        0 & \lambda_i^{-1} I_m \\
                      \end{array}
                    \right)$.

Consider the system of the form
\begin{equation}\label{linm411}
\begin{aligned}
\dot{\tilde{x}}_{ci}(t) =& \left(
               \begin{array}{cc}
                  A & 0 \\
                  G_2C & G_1 \\
               \end{array}
             \right) \tilde{x}_{ci}(t)
                            +
                             \lambda_i \left(
                               \begin{array}{c}
                                 B \\
                                 0 \\
                               \end{array}
                             \right)\tilde{u}_i(t-\tau)   ,\  i=1,\dots,N
\end{aligned}
\end{equation}
where $\tilde{x}_{ci} \in R^{(n+n_z)}$ and $\tilde{u}_i \in R^m$.

By Lemma \ref{LemAB}, there exists a matrix $\tilde{K}=(K_x,K_z)$,
where $K_x \in R^{m \times n}$ and $K_z \in R^{m \times n_z}$ such that the state feedback control law $\tilde{u}_i(t)= \tilde{K} \tilde{x}_{ci}(t)  $,  $i=1,\dots,N,$
asymptotically stabilize the system \eqref{linm411}.

Since
\begin{equation*}\label{}
\begin{split}
&\left(
  \begin{array}{cc}
    B & B \\
    0_{n_z \times m} & 0_{n_z \times m} \\
  \end{array}
\right)\left(
         \begin{array}{cc}
           K_x & 0 \\
           0 & K_z \\
         \end{array}
       \right)  = \left(
                   \begin{array}{c}
                     B \\
                     0_{n_z \times m} \\
                   \end{array}
                 \right) \left(
                           \begin{array}{cc}
                             K_x & K_z \\
                           \end{array}
                         \right) \\
\end{split}
\end{equation*}
we have $\bar{T}_{4i}  \hat{u}_{ci} (t) =  \left(
                                          \begin{array}{cc}
                                            K_x & 0 \\
                                            0 & K_z \\
                                          \end{array}
                                        \right)
 \tilde{x}_{ci}(t) $. Thus,
\begin{equation*}\label{}
\begin{split}
 \hat{u}_{ci} (t) =& \bar{T}_{4i}^{-1}\left(
                                  \begin{array}{cc}
                                    K_x & 0 \\
                                    0 & K_z \\
                                  \end{array}
                                \right)\tilde{x}_{ci}(t)\\
=&   \bar{T}_{4i}^{-1}\left(
                                  \begin{array}{cc}
                                    K_x & 0 \\
                                    0 & K_z \\
                                  \end{array}
                                \right) T_{4i} \hat{x}_{ci}(t)\\
=&\left(
    \begin{array}{cc}
      I_m & 0 \\
      0 & \lambda_i I_m \\
    \end{array}
  \right)
  \left(
    \begin{array}{cc}
       K_x & 0 \\
       0 & K_z\\
    \end{array}
  \right)
  \left(
    \begin{array}{cc}
      I_n & 0 \\
      0 & \lambda_i^{-1}I_{n_z} \\
    \end{array}
  \right)
  \hat{x}_{ci}(t)\\
=&\left(
    \begin{array}{cc}
     K_x & 0 \\
      0 & K_z \\
    \end{array}
  \right)\hat{x}_{ci}(t).\\
\end{split}
\end{equation*}
Let $K=\left(
         \begin{array}{cc}
           K_x & 0  \\
           0 & K_z \\
         \end{array}
       \right)
$. Then, we have $\hat{u}(t)= \left( I_N \otimes K \right) \hat{x}_c(t)$. Furthermore, since $\hat{u}_{c}=\bar{T}_3 \bar{T}_2 u_{c}$, $\hat{x}_c=T_3\bar{x}_c$ and $\bar{x}_c(t)=T_2x_c(t)$, we have
\begin{equation*}\label{}
\begin{split}
  u_{c} (t)  =& \bar{T}_2^{-1} \bar{T}_3^{-1} \hat{u}_{c}\\
=& \bar{T}_2^{-1} \bar{T}_3^{-1} \left( I_N \otimes K \right) \hat{x}_c(t)\\
=& \bar{T}_2^{-1} \bar{T}_3^{-1} \left( I_N \otimes K \right) T_3 T_2 x_c(t)\\
=& \bar{T}_2^{-1}  \left(
                     \begin{array}{cc}
                       I_N\otimes K_x & 0 \\
                       0 & I_N\otimes K_z \\
                     \end{array}
                   \right)
 T_2 x_c(t)\\
=& \left(
     \begin{array}{cc}
       T_1^{-1} \otimes I_m & 0 \\
       0 & T_1^{-1} \otimes I_m \\
     \end{array}
   \right)
   \left(
     \begin{array}{cc}
       I_N\otimes K_x & 0 \\
       0 & I_N\otimes K_z \\
     \end{array}
   \right)
   \left(
     \begin{array}{cc}
       T_1 \otimes I_n & 0 \\
       0 & T_1 \otimes I_{n_z} \\
     \end{array}
   \right)x_c(t)\\
=&   \left(
     \begin{array}{cc}
       I_N\otimes K_x & 0 \\
       0 & I_N\otimes K_z \\
     \end{array}
   \right)x_c(t)\\
\end{split}
\end{equation*}
The proof of the if part is then completed.

(Only if Part:)
Suppose the digraph $\mathcal{\bar{G}}$ does not satisfy Assumption \ref{Ass1.5}. Then, by Lemma 1 of \cite{Su1}
, $H$ has at least one eigenvalue at the origin. Without loss of generality, we assume
that $\lambda_l=0$. Then, by \eqref{linm3}
\begin{equation}\label{linm5}
\begin{aligned}
\dot{\hat{x}}_{cl}(t) =& \left(
               \begin{array}{cc}
                  A & 0 \\
                  0 & G_1 \\
               \end{array}
             \right) \hat{x}_{cl}(t)  + \left(
                              \begin{array}{cc}
                              0 & B  \\
                              0 & 0 \\
                              \end{array}
                            \right) \hat{u}_{cl}(t-\tau) \\
\end{aligned}
\end{equation}
Since the eigenvalues of $G_1$ coincide with those of $S$, under Assumption \ref{Ass1.1}, the system \eqref{linm5} and hence the system \eqref{linm1} cannot be asymptotically stable
regardless of the choice of $K$. The proof is thus completed.
\end{Proof}

Now, we are ready to present our result under the state feedback control law.

\begin{Theorem}\label{Theo1}
Under Assumptions \ref{Ass1.0} to \ref{Ass1.2}, \ref{Ass1.4} and \ref{Ass1.21}, there exist matrices $K_x \in R^{m \times n}$, $K_z \in R^{m \times n_z}$ such that the cooperative robust output regulation problem is solved by the distributed dynamic state
feedback control law \eqref{ctr1} with $(G_1,G_2)$ being the minimal $p$-copy internal model of $S$ if and only if Assumption \ref{Ass1.5} is satisfied.
\end{Theorem}

\begin{Proof}
Performing the coordinate transformation $\bar{z}_i(t-\tau_{com})=z_i(t)$, the state feedback control law \eqref{ctr1} becomes as follows:
\begin{equation}\label{ctr1tr}
\begin{split}
u_i(t)=&  K_x \eta_i(t) + K_z \bar{z}_i(t-\tau_{com}) \\
\dot{\bar{z}}_i(t) =&  {G}_1 \bar{z}_i(t) + {G}_2 e_{vi}(t) ,\ i=1,\dots,N.
\end{split}
\end{equation}

Then, under the state feedback control law \eqref{ctr1tr}, the undisturbed nominal closed-loop system is in the following form:
\begin{equation}\label{eqthe1}
\begin{aligned}
\dot{x}_c(t) =& \left(
               \begin{array}{cc}
                  I_N\otimes A  & 0 \\
                  H\otimes G_2C  &  I_N\otimes G_1  \\
               \end{array}
             \right) x_c(t)    +
             \left(
               \begin{array}{cc}
                H\otimes B & I_N\otimes B \\
                 0_{Nn_z \times Nm} & 0_{Nn_z \times Nm} \\
               \end{array}
             \right) \\ & \times
             \left(
               \begin{array}{cc}
                 I_N \otimes K_x & 0_{Nm \times Nn_z} \\
                 0_{Nm \times Nn} & I_N \otimes K_z \\
               \end{array}
             \right)x_c(t-\tau)
\end{aligned}
\end{equation}
where $x_c=\mbox{col}(x,\bar{z})$ with $x=\mbox{col}(x_1,\dots,x_N)$, $\bar{z}=\mbox{col}(\bar{z}_1,\dots,\bar{z}_N)$ and $K_x \in R^{m\times n}$, $K_z \in R^{m\times n_z}$.

By Lemma \ref{Lem1sta}, there exist matrices $K_x \in R^{m \times n}$ and $K_z \in R^{m \times n_z}$, such that
system \eqref{eqthe1} is asymptotically stable. The proof is thus completed by invoking Lemma \ref{Remsta}.
\end{Proof}

To study the output feedback case, we need the following lemma.
\begin{Lemma}\label{Lem2sta}
Consider the system of the form
\begin{equation}\label{Thee1}
\begin{aligned}
\dot{x}_c(t) =& \left(
               \begin{array}{ccc}
                  I_N\otimes A  & 0 & 0 \\
                  H\otimes G_2 C  & I_N \otimes G_1 & 0 \\
                  H\otimes LC  &  0  & I_N\otimes A - H\otimes LC\\
               \end{array}
             \right)   x_c(t)  + \left(
                              \begin{array}{ccc}
                                0&  I_N\otimes B  &  H\otimes B  \\
                                0& 0 & 0\\
                                0 & I_N\otimes B & H\otimes B \\
                              \end{array}
                            \right) \\& \times u_c(t-\tau)\\
x_c (\theta) =&  x_{c0}(\theta), \ \theta \in [-\tau,0]
\end{aligned}
\end{equation}
where $x_c \in R^{N(2n+n_z)}$, $u_c \in R^{3Nm}$, $(G_1,G_2)$ is the minimal p-copy internal model of $S$ as defined in \eqref{inter}, and
$x_{c0} \in \mathcal{C}\big([-\tau,0]$ $,R^{N(2n+n_z)}\big)$. Then, under Assumptions \ref{Ass1.1}, \ref{Ass1.2}, \ref{Ass1.3}, \ref{Ass1.4} and \ref{Ass1.21}, there exist matrices $K_1 \in R^{m \times n_z}$, $K_2 \in R^{m \times n}$ and $L \in R^{n \times p}$, such that under the state feedback control law $u_c(t)=Kx(t)$, where $K=\left(
                                                                                  \begin{array}{ccc}
                                                                                    I_N\otimes K_2 & 0 & 0 \\
                                                                                    0 & I_N\otimes K_1 & 0 \\
                                                                                    0 & 0 & I_N\otimes K_2 \\
                                                                                  \end{array}
                                                                                \right)$, system \eqref{Thee1} is asymptotically stable if and only if Assumption \ref{Ass1.5} is satisfied.
\end{Lemma}
\begin{Proof}
Let
$T=\left(
     \begin{array}{ccc}
       I_{Nn} & 0 & 0 \\
       0 & I_{Nn_z} & 0 \\
       -I_{Nn} & 0 & I_{Nn} \\
     \end{array}
   \right)
$ and $\bar{x}_c=Tx_c$. Then, the system \eqref{Thee1} becomes
\begin{equation}\label{tThee2}
\begin{aligned}
\dot{\bar{x}}_c(t)\!\! &=\!\! \left(\!\!\!
               \begin{array}{ccc}
                  I_N\otimes A \!\! & 0\!\! & 0 \\
                  H\otimes G_2 C \!\! & I_N \otimes G_1 \!\!& 0 \\
                  0\!\! &  0  & \!\! I_N\otimes A - H\otimes LC\\
               \end{array}
         \! \!\!   \right)\!\! \bar{x}_c(t) \!\!  + \!\! \left(\!\!\!
                              \begin{array}{ccc}
                                 H\otimes B \!\!&  I_N\otimes B  \!\!& H\otimes B   \\
                                0 \!\!& 0\!\!& 0\\
                                0\!\!& 0 \!\!& 0\\
                              \end{array}
                          \! \!\! \right)  \!\! \bar{T}u_c(t-\tau)
\end{aligned}
\end{equation}
where
$\bar{T}=\left(
                                         \begin{array}{ccc}
                                           I_{Nm} & 0 & 0 \\
                                           0 & I_{Nm}& 0  \\
                                           -I_{Nm} & 0 & I_{Nm} \\
                                         \end{array}
                                       \right)
$.

Denote $\bar{x}_c=\mbox{col}(\bar{x}_{c1},\bar{x}_{c2})$ with $\bar{x}_{c1} \in R^{N(n+n_z)}$ and $\bar{x}_{c2} \in R^{Nn}$.
Then, by Lemma \ref{Lem1sta}, under Assumptions \ref{Ass1.1}, \ref{Ass1.2}, \ref{Ass1.4} and \ref{Ass1.21}, there exist matrices $K_1 \in R^{m \times n_z}$ and $K_2 \in R^{m \times n}$, such that under the state feedback control law $\bar{u}_{c1}(t)=\hat{K}\bar{x}_{c1}(t)$, where
$\hat{K} =\left(
  \begin{array}{cc}
    I_N \otimes K_2 & 0 \\
    0 & I_N \otimes K_1 \\
  \end{array}
\right)$, the following system
\begin{equation}\label{Thee2}
\begin{aligned}
\dot{\bar{x}}_{c1}(t) =& \left(
               \begin{array}{cc}
                  I_N\otimes A  & 0  \\
                  H\otimes G_2 C  & I_N \otimes G_1 \\
               \end{array}
             \right) \bar{x}_{c1}(t)   + \left(
                              \begin{array}{cc}
                                 H\otimes B&  I_N\otimes B\\
                                0 & 0\\
                              \end{array}
                            \right) \bar{u}_{c1}(t-\tau)\\
\end{aligned}
\end{equation}
is asymptotically stable if and only if the digraph satisfies Assumption \ref{Ass1.5}. Thus, the only if part has been proved.
To show the if part, let $\bar{K}=\left(
                                                \begin{array}{ccc}
                                                  I_N\otimes K_2 & 0 & 0 \\
                                                  0 & I_N\otimes K_1 & 0 \\
                                                  0 & 0 & I_N\otimes K_2 \\
                                                \end{array} \right)$.
Then, under the state feedback control law $ u_c(t)= \bar{T}^{-1} \bar{K}\bar{x}_c(t)$, the closed-loop system of \eqref{tThee2} is as follows:
\begin{equation}\label{Eqlemc1}
\begin{split}
\dot{\bar{x}}_{c1}(t) =& \left(
               \begin{array}{cc}
                  I_N\otimes A  & 0  \\
                  H\otimes G_2 C  & I_N \otimes G_1 \\
               \end{array}
             \right) \bar{x}_{c1}(t)    + \left(
                              \begin{array}{cc}
                                 H\otimes B&  I_N\otimes B\\
                                0 & 0\\
                              \end{array}
                            \right)  \\& \times \left(
                                       \begin{array}{cc}
                                         I_N \otimes  K_2& 0 \\
                                         0 &  I_N \otimes  K_1 \\
                                       \end{array}
                                     \right)
                            \bar{x}_{c1}(t-\tau)  + \left(
                                                         \begin{array}{c}
                                                           H\otimes B \\
                                                           0 \\
                                                         \end{array}
                                                       \right) (I_N \otimes  K_2) \bar{x}_{c2}(t-\tau)\\
\dot{\bar{x}}_{c2}(t) =& (I_N\otimes A - H\otimes LC) \bar{x}_{c2}(t).
\end{split}
\end{equation}
where  $ \bar{x}_{c} = \mbox{col}(\bar{x}_{c1}, \bar{x}_{c2})$. We first note, from the proof of Theorem 2 of \cite{Su3}, that, under Assumption \ref{Ass1.3},  there exists a matrix $L$ such that the matrix $\left(I_N\otimes A - H\otimes LC\right)$ is Hurwitz. Moreover by Lemma \ref{Lem1sta}, the $ \bar{x}_{c1}$ subsystem with $\bar{x}_{c2}$ setting to zero is asymptotically stable.
Thus, by Lemma \ref{Lem2}, system \eqref{Eqlemc1} is asymptotically stable.
Furthermore, since
$\bar{x}_c(t)=Tx_c(t)$, we have
\begin{equation}\label{}
\begin{split}
&u_c(t)\\=&\bar{T}^{-1}\bar{K}\bar{x}_c(t) \\=& \bar{T}^{-1}\bar{K} Tx_c(t)\\
=& \left(
     \begin{array}{ccc}
       I_{Nm} & 0 & 0 \\
       0 & I_{Nm} & 0 \\
       I_{Nm} & 0 & I_{Nm} \\
     \end{array}
   \right) \left(
            \begin{array}{ccc}
              I_N\otimes K_2 & 0 & 0 \\
              0 & I_N\otimes K_1 & 0 \\
               0 & 0 & I_N\otimes K_2 \\
            \end{array}
          \right)
          \left(
            \begin{array}{ccc}
              I_{Nn} & 0 & 0 \\
              0 & I_{Nn_z} & 0 \\
              -I_{Nn} &  0 & I_{Nn}\\
            \end{array}
          \right) x_c(t)\\
=&\left(
     \begin{array}{ccc}
       I_N\otimes K_2 & 0 & 0 \\
       0 & I_N\otimes K_1 & 0 \\
       0 & 0 & I_N\otimes K_2 \\
     \end{array}
   \right)x_c(t) \\
   =& K x_c(t)
\end{split}
\end{equation}
 The proof is thus completed.
\end{Proof}

\begin{Theorem}\label{Theo2}
Under Assumptions \ref{Ass1.0} to \ref{Ass1.4} and \ref{Ass1.21}, there exist matrices $K_1 \in R^{m \times n_z}$, $K_2 \in R^{m \times n}$ and $L \in R^{n \times p}$ such that the cooperative robust output regulation problem is solved by the distributed dynamic output
feedback control law \eqref{ctr2} with $(G_1,G_2)$ being the minimal $p$-copy internal model of $S$ if and only if Assumption \ref{Ass1.5} is satisfied.
\end{Theorem}

\begin{Proof}
By introducing the coordinate transformation $\bar{\xi}_i(t-\tau_{com})=\xi_i(t),\ \bar{z}_i(t-\tau_{com})=z_i(t)$, the distributed dynamic output feedback control law \eqref{ctr2} becomes the
following form:
\begin{equation}\label{ctr22}
\begin{split}
u_i(t)=& K_1 \bar{z}_i(t-\tau_{com})  +  K_2 \bar{\eta}_i(t-\tau_{com}),\  i=1,\dots,N\\
\dot{\bar{z}}_i(t) =&  {G}_1 \bar{z}_i(t) + {G}_2 e_{vi}(t) \\
\dot{\bar{\xi}}_i(t) =& A \bar{\xi}_i(t) + Bu_i(t-\tau_{con}) -LC \bar{\eta}_i(t) + Le_{vi}(t)\\
\end{split}
\end{equation}
where $\bar{\eta}_i(t)=\sum_{j\in\mathcal{{N}}_i}a_{ij}\left(\bar{\xi}_i(t)-\bar{\xi}_j(t)\right)$.

Then, under the output feedback control law \eqref{ctr22}, the undisturbed nominal closed-loop system is in the following form:
\begin{equation}\label{eqqThee1}
\begin{aligned}
\dot{x}_c(t) =& \left(
               \begin{array}{ccc}
                  I_N\otimes A  & 0 & 0 \\
                  H\otimes G_2 C  & I_N \otimes G_1 & 0 \\
                  H\otimes LC  &  0  & I_N\otimes A - H\otimes LC\\
               \end{array}
             \right)   x_c(t)  + \left(
                              \begin{array}{ccc}
                                0&  I_N\otimes B  & H\otimes B   \\
                                0& 0 & 0\\
                                0 & I_N\otimes B & H\otimes B \\
                              \end{array}
                            \right) \\& \times \left(
                                      \begin{array}{ccc}
                                        I_N \otimes K_2 & 0 & 0 \\
                                        0 & I_N \otimes K_1 & 0 \\
                                        0 & 0 & I_N \otimes K_2  \\
                                      \end{array}
                                    \right)
                            x_c(t-\tau)\\
\end{aligned}
\end{equation}
where $x_c=\mbox{col}(x,\bar{z},\bar{\xi})$ with $x=\mbox{col}(x_1,\dots,x_N)$, $\bar{z}=\mbox{col}(\bar{z}_1,\dots,\bar{z}_N)$ and $\bar{\xi}=\mbox{col}(\bar{\xi}_1,\dots,\bar{\xi}_N)$.

By Lemma \ref{Lem2sta}, there exist matrices $K_1 \in R^{m \times n_z}$, $K_2 \in R^{m \times n}$ and $L \in R^{n \times p}$, such that system \eqref{eqqThee1} is asymptotically stable. The proof is thus completed by noting Lemma \ref{Remsta}.
\end{Proof}


\begin{Remark}
It is known that the cooperative output regulation problem includes the leader-following consensus problem as a special case \cite{Su1,Su3}. By the same token,
the results of this paper lead to the solution of the the leader-following consensus problem of multi-agent systems with time-delay as special cases. It is noted that, in
 \cite{HH} and \cite{Cheng},  the leader-following consensus problem of double integrator multi-agent systems with time-varying communication delays were studied under both fixed and switching communication topology. However, the control laws proposed in \cite{HH} and \cite{Cheng} need to use the speed information of the leader. Additionally, our results allow the plant to be uncertain,
 the dynamics of the leader to be different from the followers', and can reject the external disturbances.
\end{Remark}

\section{Example}\label{se}

\begin{figure}[h!]
\begin{center}
\scalebox{0.7}{\includegraphics{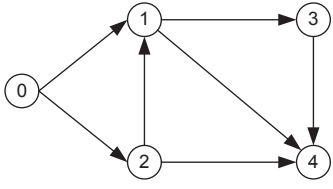}}\caption{The network topology $\mathcal{\bar{G}}$}\label{net}
\end{center}
\end{figure}

In this section, we will  illustrate our approach using the following uncertain system with input time-delay:
\begin{equation}\label{Exasys}
\begin{split}
\dot{x}_i(t) =& \left(
      \begin{array}{cc}
        0 & 1+w_{i1} \\
        0 & 0 \\
      \end{array}
    \right) x_i(t)
    + \left(
    \begin{array}{c}
    w_{i2} \\
       1 \\
        \end{array}
    \right) u_i(t-\tau_{con})
      +\left(
    \begin{array}{cc}
      0 & 0 \\
      0 & i \\
    \end{array}
  \right)v(t), \ t\geq 0,\\
e_i(t)  =& \left(
     \begin{array}{cc}
       1 & 0 \\
     \end{array}
   \right) x_i(t) - \left(
     \begin{array}{cc}
       1 & 0 \\
     \end{array}
   \right) v(t),i=1,2,3,4,\\
\end{split}
\end{equation}
with the exosystem as follows:
\begin{equation}
\begin{split}
\dot{v}(t) &= S v(t)=\left(
     \begin{array}{cc}
       0 & 1 \\
       0 & 0 \\
     \end{array}
   \right) v(t).
\end{split}
\end{equation}
The nominal system matrices are
$A=\left(
      \begin{array}{cc}
        0 & 1\\
        0 & 0 \\
      \end{array}
    \right)$,
$B=\left(
     \begin{array}{c}
       0 \\
       1 \\
     \end{array}
   \right)
$,
$E_i=\left(
    \begin{array}{cc}
      0 & 0 \\
      0 & i  \\
    \end{array}
  \right),i=1,2,3,4$,
$C= \left(
     \begin{array}{cc}
       1 & 0 \\
     \end{array}
   \right)$,
$F=\left(
     \begin{array}{cc}
       -1 & 0 \\
     \end{array}
   \right)$. The input delay $\tau_{con}=0.5s$. Here, $x_{i1}$ and $x_{i2}$ can be viewed as the position and velocity of the $i^{th}$ agent respectively and $e_{i}$ can be viewed as the tracking error of the position of the $i^{th}$ agent. The communication network topology is described in Figure \ref{net}. The matrix $H$ associated with digraph $\mathcal{\bar{G}}$ is
$H=\left(
     \begin{array}{cccc}
       2 & -1 & 0 & 0 \\
       0 & 1 & 0 & 0 \\
       -1 & 0 & 1 & 0 \\
       -1 & -1 & -1 & 3 \\
     \end{array}
   \right)
$ and the eigenvalues of $H$ are $\{3,1,2,1\}$.

It is easy to verify that Assumptions \ref{Ass1.0} to \ref{Ass1.21} are satisfied. Therefore, by Theorem \ref{Theo1} and \ref{Theo2}, the cooperative robust output regulation problem for this example can be solved by the distributed controllers of the form \eqref{ctr1} and \eqref{ctr2}, respectively.

\noindent \textbf{$(1)$ Distributed dynamic state feedback control}

The distributed dynamic state feedback controller is given as
\begin{equation}\label{ }
\begin{split}
u_i(t)=&  K_x \eta_i(t) + K_z z_i(t),\\
\dot{z}_i(t) =&  {G}_1 z_i(t) + {G}_2 e_{vi}(t-\tau_{com}) ,\ i=1,\dots,N,
\end{split}
\end{equation}
with
\begin{equation}\label{gg1}
G_1=\left(
                                                                                                                             \begin{array}{cc}
                                                                                                                               0 & 1  \\
                                                                                                                               0 & 0 \\
                                                                                                                             \end{array}
                                                                                                                           \right)\ \text{and}\
G_2=\left(
             \begin{array}{cc}
               0 \\
               1 \\
             \end{array}
           \right).
\end{equation}
Assume the communication delay $\tau_{com}=0.5s$.

Denote $ {A}_c=\left(
                      \begin{array}{cc}
                        A & 0 \\
                        G_2C & G_1 \\
                      \end{array}
                    \right)$ and $ {B}_c=\left(
                                           \begin{array}{c}
                                             B \\
                                             0 \\
                                           \end{array}
                                         \right)
                    $. By Lemma \ref{LemAB}, the desirable feedback gain is 
\begin{equation}\label{gg2}
K=(K_x,K_z)=-\nu_1^{-1} B_c^TPe^{A_c\tau}
\end{equation}
where $\tau=1$ and $P$ is the positive definite solution of the parametric ARE
\begin{equation}\label{}
A_c^T P + PA_c - PB_cB_c^T P=-\gamma P
\end{equation}
where $\gamma$ is some sufficiently small positive number.

\begin{figure}[h!]
\begin{center}
\scalebox{0.5}{\includegraphics {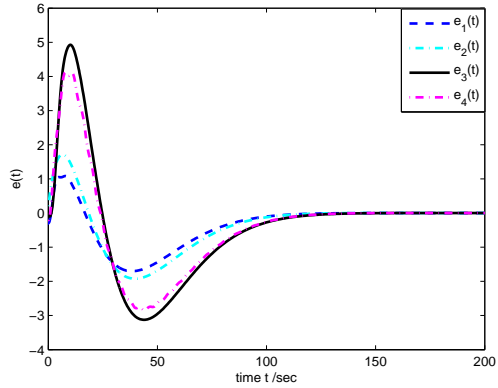}}\caption{The tracking error $e(t)$ under distributed dynamic state feedback control}\label{error}
\end{center}
\end{figure}

Figure \ref{error} shows the tracking error $e(t)$ tends to zero asymptotically where the system uncertainties are
$w=(0.05,0.03,0.05,0.01,0.02,0.08,0.05,0.04)^T$, $\nu_1=1$ and $\gamma=0.1$.

\noindent \textbf{$(2)$ Distributed dynamic output feedback control}

The distributed dynamic output feedback control law is given as
\begin{equation}\label{}
\begin{split}
u_i(t)=&  K_1 z_i(t) + K_2 \hat{\eta}_i(t) , \ i=1,\dots,N\\
\dot{z}_i(t) =&  {G}_1 z_i(t) + {G}_2 e_{vi}(t-\tau_{com}) \\
\dot{\xi}_i(t) =& A \xi_i(t) + Bu_i(t-\tau)  -LC \hat{\eta}_i(t) + Le_{vi}(t-\tau_{com}) \\
\end{split}
\end{equation}
with $(G_1,G_2)$, $(K_1,K_2)=(K_z,K_x)$ defined in \eqref{gg1} and \eqref{gg2}, respectively.
Let $\nu_2=\frac{1}{3}$, we have $L^T = \nu_2^{-1}C \tilde{P}$, where $\tilde{P}$ is the solution of the Riccati Equation
\begin{equation}\label{}
A \tilde{P} + \tilde{P} A^T + I_n - \tilde{P}C^TC\tilde{P}=0.
\end{equation}

\begin{figure}[!htb]
\begin{center}
\scalebox{0.5}{\includegraphics{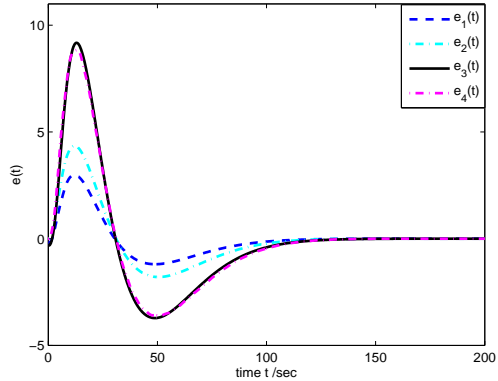}}\caption{The tracking error $e(t)$ under distributed dynamic output feedback control}\label{error2}
\end{center}
\end{figure}

Choosing $\gamma=0.1$, Figure \ref{error2} shows that the distributed dynamic output feedback controller solves the linear robust cooperative output regulation problem successfully.

To close this section, we note that this example cannot be handled by any existing methods.

\vspace{-6pt}

\section{Conclusion}\label{sc}
\vspace{-2pt}

In this paper, we have studied the cooperative robust output regulation problem of linear multi-agent
systems by the distributed internal model approach, which includes the leader-following consensus problem as a special case.
A distinguished advantage of the distributed internal model approach over the distributed observer approach in \cite{cns} is that it allows the plant parameters to be uncertain. To our knowledge, this is the first paper to handle the consensus
problem for linear uncertain multi-agent
systems with both the input delay and communication delay. Our approach can also be extended to the systems containing multiple input time-delays and state time-delays.

\end{document}